\documentclass{elsarticle}
\usepackage[a4paper, total={6in, 10in}]{geometry}
\usepackage[utf8]{inputenc}
\usepackage[center]{caption}
\usepackage{indentfirst}
\usepackage{enumerate}
\usepackage{graphicx}
\usepackage{subcaption}
\usepackage{amssymb}
\usepackage{amsmath}
\usepackage{amsthm}
\newtheorem*{theorem}{Theorem}
\newtheorem*{lemma}{Lemma}

\makeatletter
\def\ps@pprintTitle{%
  \let\@oddhead\@empty
  \let\@evenhead\@empty
  \let\@oddfoot\@empty
  \let\@evenfoot\@oddfoot
}
\makeatother

\title{A note about a transition of Ratliff and Rosenthal's order picking algorithm for rectangular warehouses}

\author[1]{Paul Revenant} \ead{paul.revenant@ens-lyon.fr}
\author[2]{Hadrien Cambazard} \ead{hadrien.cambazard@grenoble-inp.fr}
\author[2]{Nicolas Catusse} \ead{nicolas.catusse@grenoble-inp.fr}

\affiliation[1]{organization={ENS de Lyon}, country={France}}
\affiliation[2]{organization={Univ.Grenoble Alpes, CNRS, Grenoble INP, G-SCOP}, city={Grenoble}, country={France}}

\begin{document}

\begin{abstract}
    In the order picking problem, a picker has to collect a number of products in a warehouse with a minimum length tour. Ratliff and Rosenthal gave a linear algorithm for solving the order picking problem in the case where the warehouse has two cross aisles. Their algorithm allows the tour to double cross an entire aisle. We prove that, in rectangular warehouses, there always exists a minimum length tour which doesn't double cross an aisle.
\end{abstract}
\maketitle

\section{Introduction}
\begin{figure}[h!]
    \centering
    \includegraphics[scale=0.5]{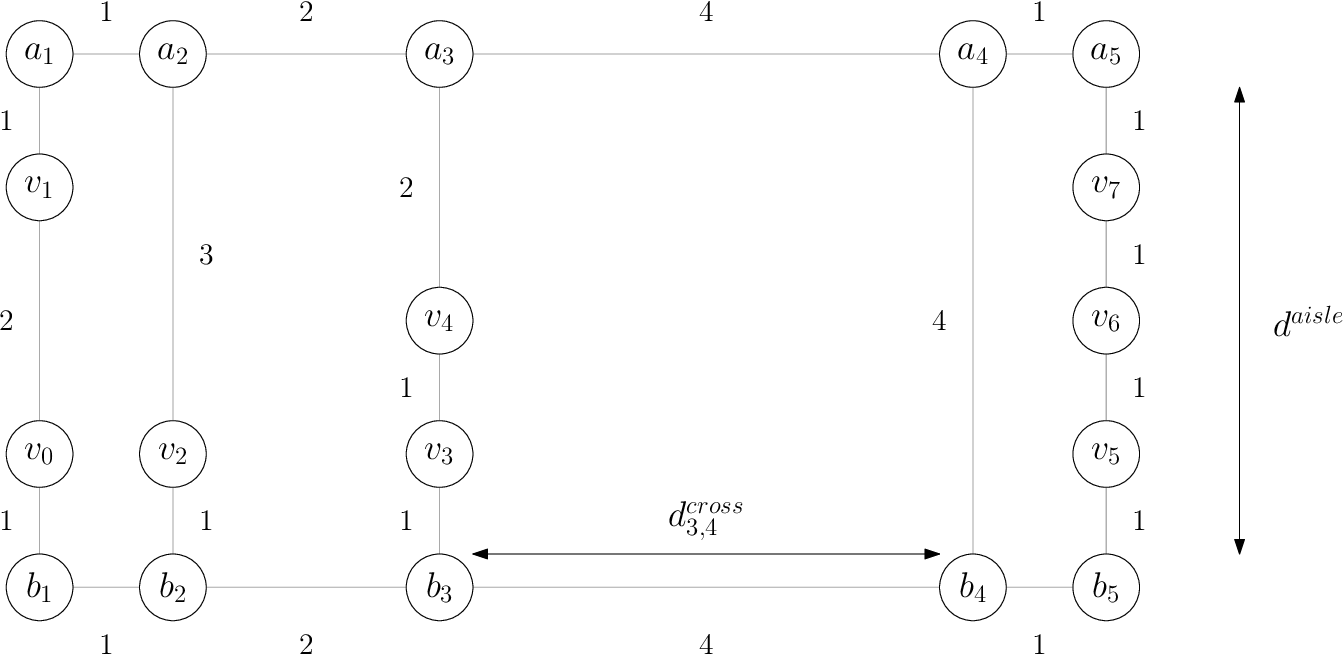}
    \caption{An example of a picking instance in a rectangular warehouse with five aisles and two cross aisles.}
    \label{fig:example_warehouse}
\end{figure}

Order picking is a problem widely studied for decades in operations research, due to its importance in supply chains. In 1983, Ratliff and Rosenthal published a linear algorithm for warehouses with two cross aisles (see \cite{ratliff} and \cite{note_lin}). The algorithm is presented for rectangular warehouses similar to the one of Figure \ref{fig:example_warehouse}. The authors also mentioned it can be applied to non-rectangular warehouses where aisles (and cross aisles) do not necessarily have the same length. They use a dynamic-programming approach, based on the possible edge configurations of the tour along the aisles and the cross aisles. In this note, we show that one of these configurations, namely double-crossing an aisle, is unnecessary in rectangular warehouses. Note that Ratliff and Rosenthal already mention this fact (without proof)  as a side remark, although what they mean is not entirely clear.

We consider a warehouse graph $G$ having two horizontal cross aisles and $n$ vertical aisles, with $n \geq 1$. The vertices of $G$ include $a_1, \ldots, a_n$ and $b_1, \ldots, b_n$, which are respectively the upper and lower intersections between the aisles and the cross aisles. In what follows, we assume that $G$ has a rectangular shape, which means that all the aisles have the same length $d^{aisle}$, and that the two cross aisles between aisles $j$ and $j+1$ have the same length $d_{j, j+1}^{cross}$, for $1 \leq j < n$. Let us denote $P = \{v_0, v_1, \ldots, v_k, \ldots, v_m\}$ the set of vertices to be visited during the tour on $G$, by identifying $v_0$ with the depot, and $\{v_1, \ldots, v_m\}$ with the products to be picked up. Except for $v_0$ which can coincide with a vertex $a_j$ or $b_j$, all the vertices of $P$ are strictly contained into an aisle of the warehouse.

\begin{figure}[h!]
    \centering
    \includegraphics[scale=0.4]{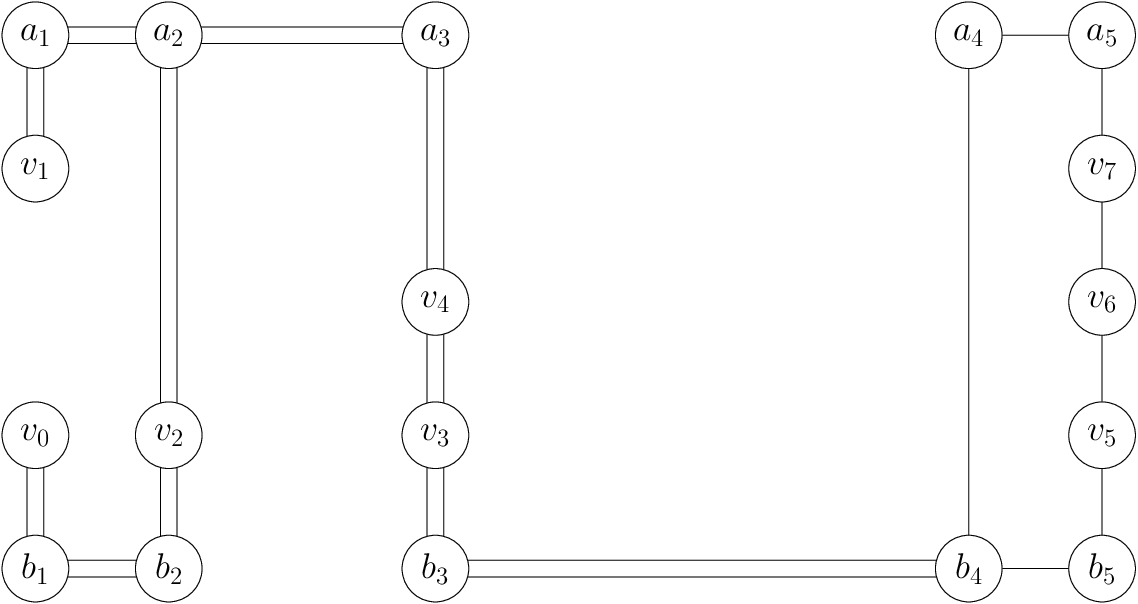}
    \caption{A non optimal tour subgraph corresponding to the previous instance,\\with two double edges crossing the aisles $(a_2, b_2)$ and $(a_3, b_3)$.}
    \label{fig:example_warehouse_non_opti}
\end{figure}

\begin{figure}[h!]
    \centering
    \includegraphics[scale=0.4]{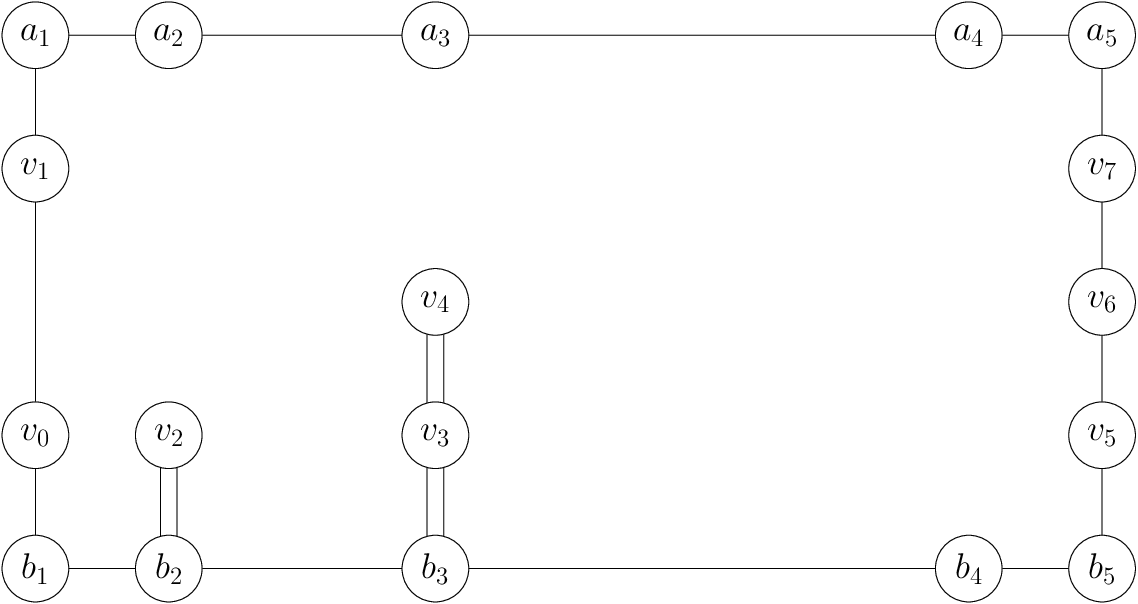}
    \caption{A minimum length tour subgraph solving the previous instance.\\It has no double edge crossing an entire aisle.}
    \label{fig:example_warehouse_opti}
\end{figure}

\begin{figure}[h!]
    \centering
    \includegraphics[scale=0.46]{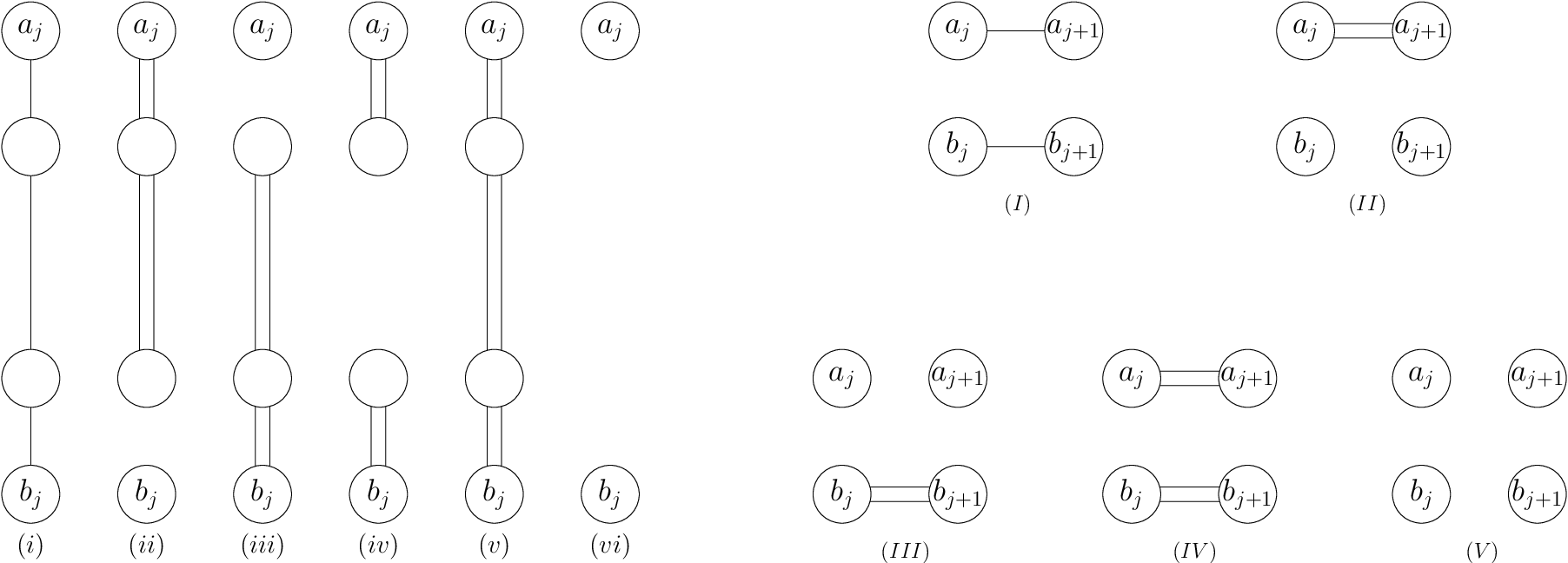}
    \caption{The possible edge configurations for vertical edges (on the left) and for\\horizontal edges (on the right). The double edge across an aisle is the configuration $(v)$. }
    \label{fig:arc_config}
\end{figure}

\begin{figure}[h!]
    \centering
    \captionsetup{justification=centering}
    \includegraphics[scale=0.45]{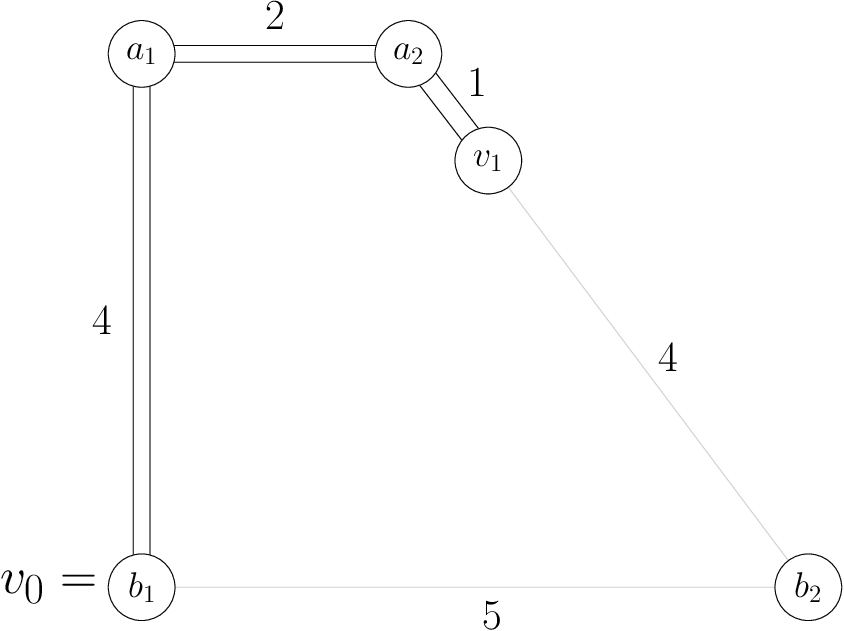}
    \caption{An example of a minimum length tour subgraph visiting $v_0$ and $v_1$ using a necessary vertical double edge crossing the aisle $(a_1, b_1)$ in a non rectangular warehouse.}
    \label{fig:mandatory_double_edge}
\end{figure}

\newpage
Ratliff and Rosenthal list all the possible edge configurations along a vertical aisle and along a horizontal cross aisle (see Fig.~\ref{fig:arc_config}). One of them seems to be unnecessary: the vertical double edge crossing an entire aisle (configuration (v) on the left of Figure \ref{fig:arc_config}). Note that in the general case, where the warehouse is not rectangular, this configuration is required, as shown in Fig.~\ref{fig:mandatory_double_edge}. But in the case of a rectangular warehouse, there always exists an optimal tour that does not require (v). We introduce some preliminary results from \cite{ratliff} in Section \ref{requisites} and state the main result in Section \ref{result}.

\section{A few prerequisites}
\label{requisites}
The proofs of the lemmas presented in this section can be found in \cite{ratliff}. Define a subgraph $T \subseteq G$ as a \textit{tour subgraph} if $(i)$ all vertices $v_0, \ldots, v_m$ have positive degree in $T$; $(ii)$ $T$ is connected; and $(iii)$ every vertex in $T$ has even degree. Given a subgraph $L \subseteq G$, a subgraph $T \subseteq L$ is a \textit{$L$ partial tour subgraph} (PTS) if there exists a subgraph $C \subseteq G - L$ such that $T \cup C$ is a tour subgraph of $G$. The subgraph $C$ is called a \textit{completion} of $T$.

\begin{figure}[h!]
    \centering
    \includegraphics[scale=0.5]{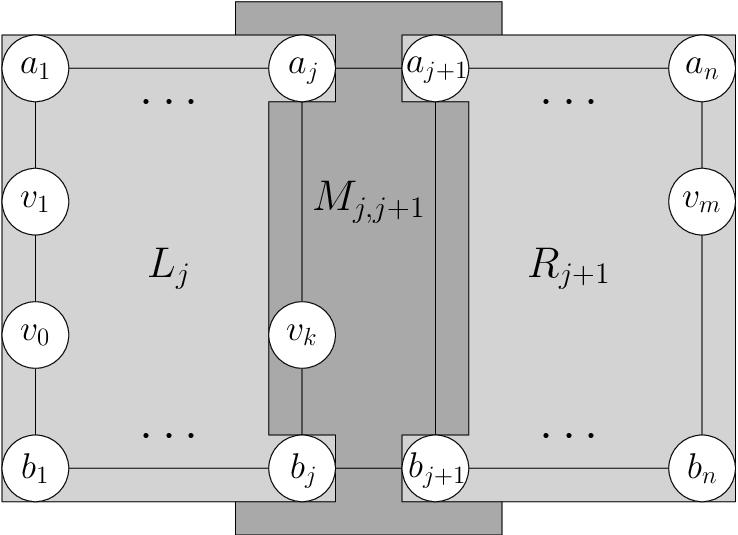}
    \caption{Illustration of the definitions of $L_j$, $R_j$ and $M_{j, j+1}$.}
    \label{fig:def_lj}
\end{figure}

Let us define the following subgraphs of $G$. For $1 \leq j \leq n$, $L_j$ (resp $R_j$) consists in the vertices $a_j$ and $b_j$ (resp. $a_{j+1}$ and $b_{j+1}$), with everything in $G$ strictly to the left (resp. to the right) of aisle $j$ (resp. $j+1$). Moreover, for $1 \leq j < n$, we denote $M_{j, j+1}$ the subgraph induced by all vertices of aisles $j$ and $j+1$. It includes vertices $a_j$, $b_j$ and everything in between, together with $a_{j+1}$, $b_{j+1}$ and everything in between (see Fig.~\ref{fig:def_lj}).

Let $S_j$ be any of the subgraphs among $L_j$ or $R_j$, for $1 \leq j \leq n$. We have the following useful characterization of a $S_j$ PTS:
\begin{lemma} The subgraph $T \subseteq S_j$ is a $S_j$ PTS if and only if 
\begin{enumerate}[$(i)$]
    \item for all $v_i \in S_j$, the degree of $v_i$ is positive in $T$;
    \item every vertex in $T$ has an even degree, except possibly for $a_j$ and $b_j$; and
    \item $T$ has either no connected component, or a single connected component containing at least one of $a_j$ and $b_j$, or two connected components with $a_j$ in one component and $b_j$ in the other.
\end{enumerate}
\end{lemma}

Two $S_j$ PTS are said to be equivalent if any completion of one PTS is also a completion of the other. It corresponds to the following characterization:
\begin{lemma} The two $S_j$ PTS $T_1$ and $T_2$ are equivalent if and only if
\begin{enumerate}[$(i)$]
    \item $a_j$ has the same degree parity (even but not zero, odd or zero) in both $T_1$ and $T_2$, and $b_j$ has the same degree parity in both $T_1$ and $T_2$; and
    \item both $T_1$ and $T_2$ have no connected component, or both have a single connected component containing at least one of $a_j$ and $b_j$, or both have two connected components with $a_j$ in one component and $b_j$ in the other.
\end{enumerate}
\end{lemma}

We denote the class equivalence of a $S_j$ PTS as (degree parity of $a_j$, degree parity of $b_j$, number of connected components), where the degree parity of a vertex can be $E$ (even but not zero), $U$ (uneven) or $0$. There exist exactly seven equivalence classes for a $S_j$ PTS, namely $(U,U,1C)$, $(E,E,1C)$, $(0,E,2C)$, $(E,0,2C)$, $(E,E,2C)$, $(0,0,0C)$, $(0,0,1C)$. In particular, we will often use the fact that if a $S_j$ PTS $T$ doesn't belong to the class $(U,U,1C)$, then $a_j$ and $b_j$ have an even degree in $T$.

In addition, we define $d(T)$ as the sum of the edge distances of a subgraph $T \subseteq G$.

\section{Discarding vertical double edges along an entire aisle}
\label{result}
In this section, we prove the following:

\begin{theorem}
There exists a minimum length tour subgraph $T \subseteq G$ such as $T$ has no vertical double edge crossing an entire aisle.
\end{theorem}

Before presenting the sketch of the proof, let us point out some preliminary remarks.

\begin{itemize}
    \item Firstly, the case where all the vertices of $P$ lie on the same aisle doesn't contradict our statement. Indeed, the optimal tour subgraph in this case is a double edge between the two extreme points of $P$. Since we imposed that the vertices of $P$ can not coincide with a vertex $a_j$ or $b_j$ (except possibly for the vertex $v_0$), this double edge doesn't cross the entire aisle. Hence, we suppose that there exist vertices of $P$ in at least two aisles.

    \item Secondly, we can assert without loss of generality that there exist vertices of $P$ on both the far left and the far right aisle (since a minimum length tour would never visit these aisles otherwise, due to the rectangular shape of the warehouse).

    \item Thirdly, if $T$ is a tour subgraph having a double edge along the aisle $j$, we can consider without loss of generality that this aisle is not at the far right of the warehouse, i.e. $j<n$. Indeed, in the case where $j=n$, we simply flip the warehouse by a vertical symmetry, in order to get a vertical double edge along the first aisle. The aim of this observation is to ensure that the aisle $j+1$ exists in $G$. 
\end{itemize}

\vskip 10pt

Considering a tour subgraph $T \subseteq G$ having at least one vertical double edge, we will make a reduction into a tour subgraph $T'$ having strictly less vertical double edge than $T$, with $T'$ not longer than $T$. We list the conditions to check when constructing such a reduction:
\begin{enumerate}[$(i)$]
    \item $T'$ is not longer than $T$, i.e. $d(T') \leq d(T)$,
    \item $T'$ has strictly less vertical double edge than $T$,
    \item $T'$ visits all vertices of $P$,
    \item $T'$ is connected,
    \item $T'$ has only vertices of even degree.
\end{enumerate}

Hence, starting with any minimum length tour subgraph and applying this reduction iteratively will eventually lead to a minimum length tour subgraph without a vertical double edge.

\vskip 10pt

\begin{proof}
Let $T \subseteq G$ a tour subgraph having a vertical double edge between $a_j$ and $b_j$. Due to the preliminary remarks, we can assume that $G$ has at least two aisles, with products to be picked up at both the far left and the far right aisle, and that the aisle $j$ is not at the far right of the warehouse, such that the aisle $j+1$ exists in $G$.

We decompose the tour subgraph into $T = T_L \cup T_M \cup T_R$, with $T_L = T \cap L_j$, $T_M = T \cap M_{j, j+1}$ and $T_R = T \cap R_{j+1}$. In the reduction presented below, we will only modify $T_M$ into $T_M' \subseteq M_{j,j+1}$, and set $T'=T_L \cup T_M' \cup T_R$. Hence, all the figures will only show the modifications performed on $T_M$. We will always construct $T_M'$ with edge configurations appearing in the set of possible edge configurations of Fig.~\ref{fig:arc_config}. Hence, the conditions listed in the sketch of the proof become:
\begin{enumerate}[$(i)$]
\item $d(T_M') \leq d(T_M)$,
\item $T_M'$ has strictly less vertical double edge than $T_M$,
\item $T_M'$ visits all the vertices of $P$ lying on $M_{j,j+1}$,
\item $T_L \cup T_M' \cup T_R$ is connected,
\item each vertex $a_j$, $b_j$, $a_{j+1}$, $b_{j+1}$ has the same degree parity (even or odd) in $T_M$ and in $T_M'$.
\end{enumerate}

Let us consider the equivalence class $c_L$ of $T_L$ (as a $L_j$ partial tour subgraph) and the equivalence class $c_R$ of $T_R$ (as a $R_{j+1}$ partial tour subgraph). Pay attention to the fact that even if $a_j$ and $b_j$ are linked by a vertical double edge in $T$, they might be in two distinct connected components in $T_L$.
The reduction will distinguish between three cases.
\begin{description}

\item[First case: $c_L = (U,U,1C)$.] \hfill \\
In this case, we note in particular that $a_j$ and $b_j$ are in the same connected component in $T_L$. It suffices to remove from $T_M$ a double edge between two consecutive vertices of aisle $j$ to obtain $T_M'$ (see Fig.~\ref{fig:reduc1}). Both the length of the tour and the number of vertical double edges have decreased, and the tour still visits all the vertices of $P$. Moreover, $T'$ is connected. Indeed, $T$ is connected and differs from $T'$ by an edge linking two vertices which remains in the same connected component. Finally, we observe that the parity of the degrees of $a_j$, $b_j$, $a_{j+1}$ and $b_{j+1}$ remains the same.

\begin{figure}[h!]
    \centering
    \includegraphics[scale=0.5]{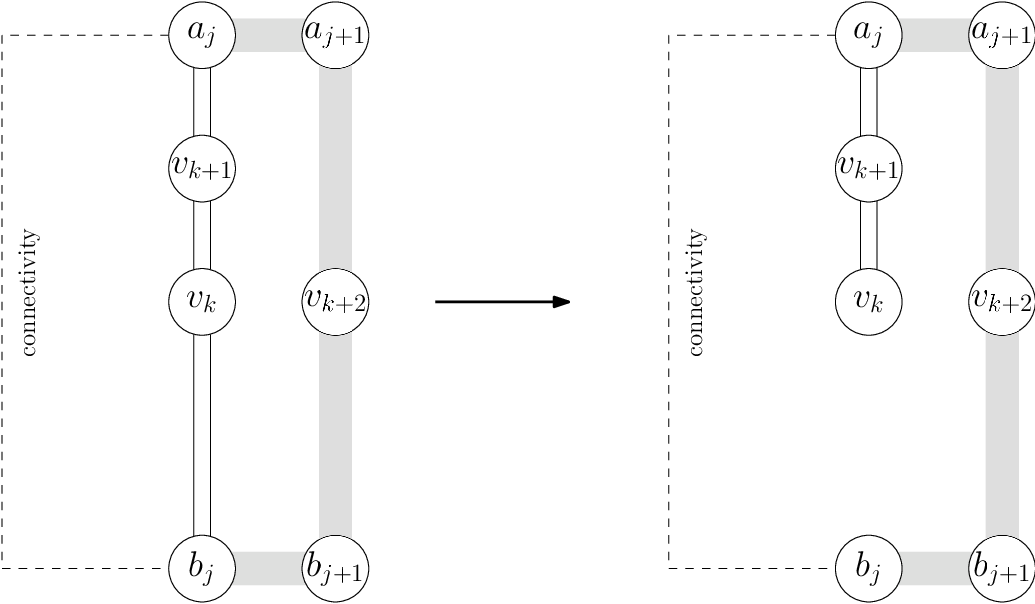}
    \caption{Reduction in the case where $c_L=(U,U,1C)$. The gray edges represent any possible edge configuration of Fig.~\ref{fig:arc_config}, they are not modified during the reduction.}
    \label{fig:reduc1}
\end{figure}

\item[Second case: $c_L \ne (U,U,1C)$ and $c_R=(U,U,1C)$.] \hfill \\
Since every vertex has an even degree in $T$, this assumption implies, in particular, that $a_j$ and $b_j$ have an even degree in $T_M$, while $a_{j+1}$ and $b_{j+1}$ have an odd degree in $T_M$. Because of the double edge between $a_j$ and $b_j$, we are left with two possibilities for the top cross aisle $(a_j, a_{j+1})$ and the bottom cross aisle $(b_j, b_{j+1})$ between aisles $j$ and $j+1$: both can have either a horizontal double edge, or no edge at all. Moreover, because $T$ is connected and has to visit the two extreme aisles of the warehouse, at least one of these cross aisles is visited. We suppose, by symmetry, that the top cross aisle is a double edge. Because of the odd parity of the degrees of $a_{j+1}$ and $b_{j+1}$ in $T_M$, the aisle $j+1$ is crossed by a single edge. It all comes to the situation depicted on the left of Fig.~\ref{fig:reduc2}. We note that $d(T_M)$ is at least $3d^{aisle}+2d^{cross}_{j, j+1}$. 

\begin{figure}[h!]
    \centering
    \includegraphics[scale=0.5]{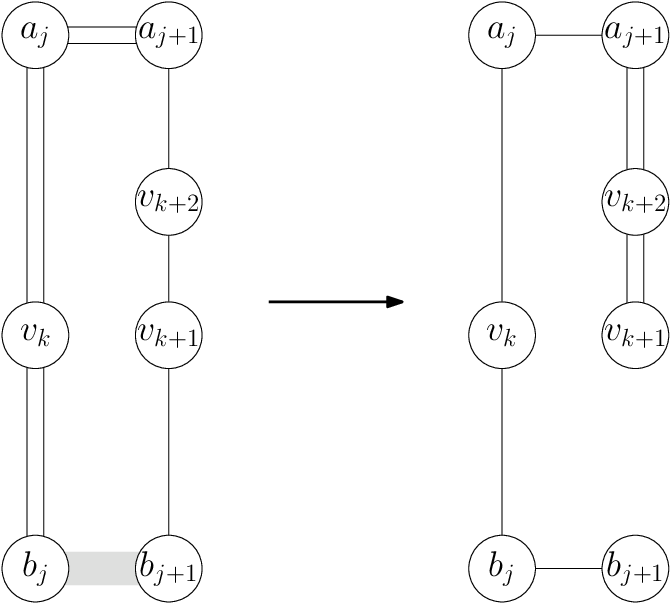}
    \caption{Reduction in the case where $c_L \ne (U,U,1C)$ and $c_R = (U,U,1C)$. The gray edge represents either a horizontal double edge or no edge at all.}
    \label{fig:reduc2}
\end{figure}

We define $T_M'$ as the subgraph of $M_{j, j+1}$ having a single edge along the aisle $j$, a single edge on both cross aisles between aisles $j$ and $j+1$, and a vertical double edge along the aisle $j+1$, except between two consecutive vertices (see Fig.~\ref{fig:reduc2}). We check that $d(T_M') < 3d^{aisle}+2d^{cross}_{j, j+1} = d(T_M)$, and that $T_M'$ doesn't have any vertical double edge, visits all vertices of $P$ in $M_{j, j+1}$, and doesn't break the connectivity of the tour. Indeed, $T_M'$ has a single connected component, which contains all the vertices of $M_{j,j+1}$. Furthermore, the degree of $a_j$, $b_j$, $a_{j+1}$, $b_{j+1}$ has the same parity in $T_M$ and in $T_M'$.

\item[Third case: $c_L \ne (U,U,1C)$ and $c_R \ne (U,U,1C)$.] \hfill \\
It implies that all the vertices $a_j$, $b_j$, $a_{j+1}$, $b_{j+1}$ have an even degree in $T_M$. Like in the previous case, it forces at least one of the two cross aisles to have a double edge, and we can suppose by symmetry that it is the case for the top cross aisle. The situation is drawn on Fig~\ref{fig:square}. We define $T_M'$ as the subgraph of $M_{j, j+1}$ having one edge along the two aisles and the two cross aisles. We check that $d(T_M') = 2d^{aisle} + 2d_{j, j+1}^{cross} \leq d(T_M)$. Note that $T_M'$ has no vertical double edge and visits all the vertices of $M_{j, j+1}$. Furthermore, the connectivity of the tour is preserved, and $a_j$, $b_j$, $a_{j+1}$, $b_{j+1}$ have an even degree in $T_M'$ and $T_M$.

\begin{figure}[h!]
    \centering
    \includegraphics[scale=0.5]{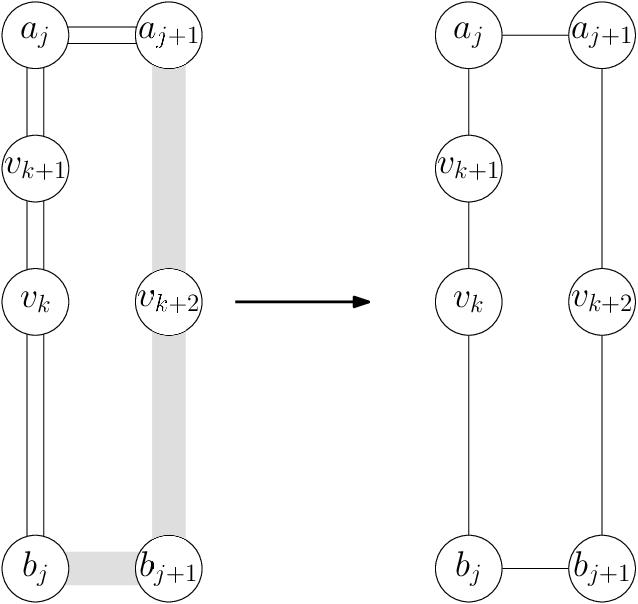}
    \caption{Reduction in the case where $c_L$ and $c_R$ are different from $(U,U,1C)$. The bottom gray edge represents a horizontal double edge or no edge at all. The vertical gray edge represents any vertical edge configuration of Fig.~\ref{fig:arc_config} respecting the even degree of $a_{j+1}$ and $b_{j+1}$ in $T_M$.}
    \label{fig:square}
\end{figure}
\end{description} 

Finally, in all three cases considered above, we have a valid reduction to strictly decrease the number of vertical double edges of the tour subgraph without increasing the length of the tour. Thus, starting with a minimum length tour subgraph and applying this reduction iteratively will eventually lead to a minimum length tour subgraph without a vertical double edge.
\end{proof}

\section{Conclusion}

\begin{figure}[h!]
    \centering
    \begin{subfigure}[b]{0.45\textwidth}
        \centering
        \includegraphics[scale=0.75]{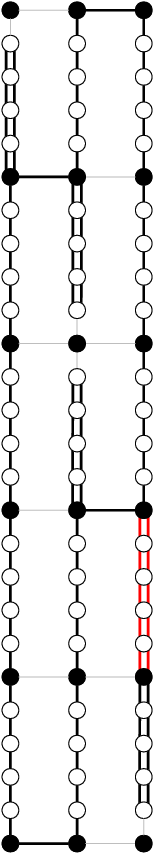}
        \caption{}
        \label{fig:two_double_edges}
    \end{subfigure}
    \hfill
    \begin{subfigure}[b]{0.45\textwidth}
        \centering
        \includegraphics[scale=0.75]{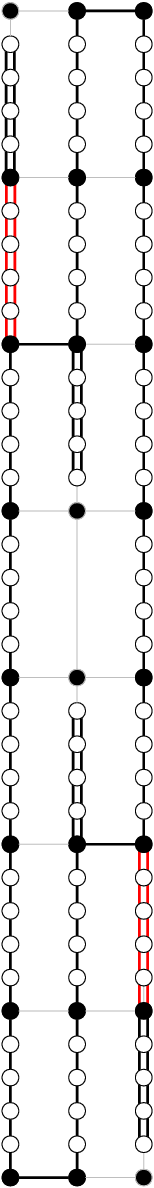}
        \caption{}
        \label{fig:second_image}
    \end{subfigure}
    \caption{(a) One necessary vertical double edge (red) in a warehouse with six cross aisles and three aisles, all filled with products (circles). The intersections are represented by solid black circles. (b) Two necessary vertical double edges (red) in a warehouse with eight cross aisles and three aisles, all filled with products except for the middle center aisle.}
    \label{fig:both_images}
\end{figure}

We have proven that considering vertical double edges along an aisle is unnecessary in rectangular warehouses with two cross aisles. A natural generalization is to consider multiple cross aisles \cite{de_koster_middle, de_koster_multiple, catusse_cambazard}. Double edges are obviously mandatory when the warehouse is made up of a single aisle and multiple cross aisles for feasibility reasons since the only tour goes up and down the warehouse. But it appears that double edges are more generally necessary in optimal solutions for rectangular warehouses with more than two cross aisles. Two examples with mandatory double edges are shown in Fig.~\ref{fig:both_images}. Nevertheless, bounding the number of double edges and restricting their locations might lead to practical improvements for the exact resolution of the order picking problem in large warehouses or the rectilinear traveling salesman problem.

\bibliographystyle{unsrt}
\bibliography{articles}

\end{document}